\documentclass[11pt]{amsart}

\usepackage[utf8]{inputenc}
\usepackage{amsmath, amssymb,amsthm,tikz}
\usepackage{tikz-cd}
\usepackage{tikz}
\usepackage{todonotes}
\usepackage{mathrsfs}
\usepackage{extarrows}
\usepackage{graphicx}
\usepackage{thmtools}
\usepackage{mathtools}
\usepackage{hyperref}
\usepackage{multirow}
\usepackage{esvect} 

\usepackage{nicematrix}
\usepackage{adjustbox}
\usepackage[margin=1.48in]{geometry}

\usepackage[english]{babel}
\usepackage{comment}

\usepackage[toc,page]{appendix}

\newcommand{\NN}{\mathbb{N}}
\newcommand{\ZZ}{\mathbb{Z}}
\newcommand{\RR}{\mathbb{R}}
\newcommand{\CC}{\mathbb{C}}
\newcommand{\QQ}{\mathbb{Q}}

\newcommand{\cC}{\mathcal{C}}
\newcommand{\cD}{\mathcal{D}}
\newcommand{\cH}{\mathcal{H}}
\newcommand{\cbT}{\mathfrak{T}}

\newcommand{\cB}{\mathcal{B}}

\newcommand{\bbG}{\mathbb{G}}

\newcommand{\bpi}{{\widetilde{\pi}}}

\newcommand{\bbbT}{{\underline{\bbT}}}

\newcommand{\bbT}{\mathbb{T}}

\newcommand{\cE}{\mathcal{E}}
\newcommand{\cF}{\mathcal{F}}
\newcommand{\cG}{\mathcal{G}}
\newcommand{\cJ}{\mathcal{J}}
\newcommand{\cL}{\mathcal{L}}
\newcommand{\cM}{\mathcal{M}}
\newcommand{\cO}{\mathcal{O}}
\newcommand{\cP}{\mathcal{P}}
\newcommand{\cS}{\mathcal{S}}
\newcommand{\cT}{\mathcal{T}}

\newcommand{\mini}{\mathrm{min}}
\newcommand{\rS}{\mathrm{S}}

\newcommand{\sign}{\operatorname{sign}}
\newcommand{\rI}{\mathrm{I}}
\newcommand{\rII}{\mathrm{II}}
\newcommand{\maxi}{\mathrm{max}}
\newcommand{\mult}{\mathrm{m}}

\newcommand{\Ind}{\operatorname{Ind}}

\newcommand{\rU}{\mathrm{U}}

\newcommand{\bG}{\mathbf{G}}

\newcommand{\an}{\operatorname{an}}

\newcommand{\tame}{\operatorname{t}}

\newcommand{\unip}{\mathrm{u}}

\newcommand{\rr}{\operatorname{r}}
\newcommand{\ii}{\operatorname{i}}

\newcommand{\cInd}{\operatorname{c-Ind}}
\newcommand{\Cent}{\operatorname{C}}

\newcommand{\Nor}{\operatorname{N}}
\newcommand{\fR}{\mathfrak{R}}
\newcommand{\fB}{\mathfrak{B}}

\newcommand{\Irr}{\operatorname{Irr}}

\newcommand{\fs}{\mathfrak{s}}

\newcommand{\fX}{{\mathfrak{X}}}

\newcommand{\der}{{\mathrm{der}}}
\newcommand{\Lie}{{\mathrm{Lie}}}
\newcommand{\rN}{{\mathrm{N}}}

\newcommand{\nr}{{\mathrm{nr}}}

\newcommand{\enh}{{\mathrm{e}}}

\newcommand{\bbbG}{{\underline{\mathbb{G}}}}

\newcommand{\rZ}{{\mathrm{Z}}}

\newcommand{\rA}{\operatorname{A}}

\newcommand{\GL}{\mathrm{GL}}
\newcommand{\GSpin}{\mathrm{GSpin}}
\newcommand{\SL}{\operatorname{SL}}

\newcommand{\rG}{\operatorname{G}}

\newcommand{\Sp}{\operatorname{Sp}}
\newcommand{\SO}{\operatorname{SO}}

\newcommand{\cusp}{\mathrm{cusp}}

\newcommand{\Sc}{\mathrm{Sc}}

\def\infl{\mathrm{infl}}

\def\bq{{\mathbf{q}}}

\def\triv{\mathrm{triv}}
\def\Unip{\mathfrak{U}}

\def\cR{\mathcal{R}}
\def\cV{\mathcal{V}}

\def\rO{{\mathrm{O}}}

\def\rA{{\mathrm{A}}}
\def\rB{{\mathrm{B}}}
\def\rC{{\mathrm{C}}}
\def\rD{{\mathrm{D}}}

\def\rG{{\mathrm{G}}}
\def\IC{{\mathrm{IC}}}
\def\Mp{{\mathrm{Mp}}}

\def\reduit{{\mathrm{r}}}
\def\ri{{\mathrm{i}}}
\def\tame{{\mathrm{t}}}
\def\scc{\mathrm{sc}}
\def\Fr{\mathrm{Fr}}
\def\an{{\mathrm{an}}}
\def\rr{{\mathrm{r}}}
\def\temp{{\mathrm{t}}}
\def\cc{{\mathrm{c}}}
\def\scc{{\mathrm{sc}}}

\def\cD{{\mathcal{D}}}
\renewcommand{\tilde}{\widetilde}

\numberwithin{equation}{subsection}
\newtheorem{theorem}[equation]{Theorem}
\newtheorem{proposition}[equation]{Proposition}

\theoremstyle{definition}
\newtheorem{definition}[equation]{Definition}
\newtheorem{remark}[equation]{Remark}
\newtheorem{example}[equation]{Example}

\title[Correspondences between affine Hecke algebras and applications]{Correspondences between affine Hecke algebras and applications}

\author{Anne-Marie Aubert}

\address{Sorbonne Universit\'e and Universit\'e Paris Cit\'e, CNRS,
IMJ-PRG, F-75005 Paris, France\\
anne-marie.aubert@imj-prg.fr}

\date{\today}

\begin{document}

\maketitle

\begin{abstract}
We review the construction of generalized affine Hecke algebras attached to Bernstein series of both smooth irreducible and enhanced $L$-parameters of $p$-adic reductive groups and apply it to the study of the Howe correspondence.

\end{abstract}

\vspace*{12pt}

\tableofcontents

\section{Introduction} \label{sec:intro}

Affine Hecke algebras are known to be very useful to describe the smooth complex representations of reductive $p$-adic groups, in terms of the supercuspidal representations
of their Levi subgroups. Their use in the classification of smooth irreducible  representations started with the seminal papers of Kazhdan and Lusztig \cite{KL} for Iwahori-spherical representations and of Lusztig \cite{LupadicI}, \cite{LupadicII} for unipotent representations.

For $G$ a $p$-adic group, the theory of the Bernstein Center of \cite{Bernstein-centre} provides a decomposition of the category $\fR(G)$ of smooth representations of $G$ as the direct product of a family of full subcategories $\mathfrak{R}^{\mathfrak{s}}(G)$,  that are called Bernstein blocks.  Here $\fs$ is the equivalence class of a pair $(L,\sigma)$, where $L$ is a Levi subgroup of $G$ and $\sigma$ is an irreducible supercuspidal representation of $L$ (see {\S}\ref{sec:Bernstein} for more details). By \cite{Heiermann}, if  $G$ is symplectic group or an orthogonal group, each Bernstein block is  Morita equivalent to the category of modules of an \textit{extended affine Hecke algebra} $\cH(G,\fs)$, that is, the crossed product of an affine Hecke algebra by the group algebra of a finite group. By \cite{Solleveld-endomorphism-algebra}, for $G$ an arbitrary $p$-adic group,  and every $\fs=[L,\sigma]_G$ such that the restriction of $\sigma$ to $L^1$  is multiplicity free, where $L^1$ is the subgroup of $L$ generated by all its compact its subgroups, the Bernstein block $\mathfrak{R}^{\mathfrak{s}}(G)$ is Morita equivalent to the category of modules of an algebra $\cH(G,\fs)$, which is  the crossed product of an affine Hecke algebra by a ($2$-cocycle) twisted group algebra of a finite group. 

In recent works, we showed that  affine Hecke algebras play a similar role on the Galois side
of the local Langlands correspondence. In \cite{AMS1}, we introduced a notion of \textit{cuspidality} for enhanced 
Langlands parameters for $G$, which conjecturally puts irreducible supercuspidal representations of $G$ in bijection with such enhanced $L$-parameters. 
We also define \textit{cuspidal support maps} and \textit{Bernstein series} for enhanced $L$-parameters, in analogy with Bernstein’s theory of representations of $p$-adic groups.
In \cite{AMS3}, to every Bernstein series $\Phi_\enh^{\fs^\vee}(G)$ of enhanced Langlands
parameters for $G$ we canonically associated an affine Hecke algebra $\cH(G^\vee,\fs^\vee)$, again possibly extended by the twisted group algebra of a finite group. We proved that the simple modules of this algebra are naturally in
bijection with the members of the Bernstein series, and that the set of central characters of
the algebra is naturally in bijection with the collection of cuspidal supports of these
enhanced Langlands parameters. These bijections send tempered or (essentially) square-integrable
representations to the expected kind of Langlands parameters.

Furthermore, in \cite{AMS3} and \cite{AMS4}, we proved that for many reductive $p$-adic groups (in particular for classical groups), if a Bernstein series $\Irr(\fR^\fs(G))$ for
$G$ corresponds to a Bernstein series $\Phi_\enh^{\fs^\vee}(G)$ of enhanced Langlands parameters for $G$ via the
local Langlands correspondence, then the algebras  $\cH(G,\fs)$ and $\cH(G^\vee,\fs^\vee)$ are Morita equivalent. 

We review the construction of these algebras, and explain how these results can be apply to describe the theta correspondence for $p$-adic groups. 
The theory of this correspondence, also known as the Howe correspondence, was  initiated by Roger Howe and builds upon the work of Weil in a representation theoretic treatment: a key role in the theory is played by the Weil representation (so-called, oscillator representation), a representation of a non-linear-algebraic group (the metaplectic group), described notably in \cite{AP}. It turned out that this correspondence has many important uses in local representation theory. Furthermore, its global counterpart (which agrees with local theta correspondence locally) provides one of a very few direct ways to explicitly construct automorphic forms.

We describe some of the recent results on the description of the Howe correspondence over both finite fields and non-archimedean local fields. For $(G_n,G'_{n'})$ an irreducible reductive dual pair, formed of a $p$-adic symplectic group $\Sp_{2n}(F)$ and a  $p$-adic orthogonal group $\rO_{2n'}(F)$, we show in Theorem~\ref{thm:corrAHA} that the compatibility the Howe correspondence  with parabolic induction leads to a collection of correspondences between simple modules of algebras $\mathcal{H}(G_n,\mathfrak{s})$ and $\mathcal{H}(G'_{n'},\mathfrak{s'})$, where $\fs'$ is explicitely described in term of $\fs$. 
We consider in particular the case of tempered representations and derive in Theorem~\ref{thm:Cstar} a correspondence at the level of $C^\ast$-algebras when $n'=n$ or $n'=n+1$. 

The study of the Howe correspondence has been also heavily influenced by the recent developments in the framework of the Langlands program, so that new results are formulated in this spirit (see notably \cite{AG}, \cite{GS}, and \cite{GI}). As a consequence of our results, we show that the Howe correspondence induces a collection of correspondences between simple modules of the algebras $\cH(G^\vee_n,\fs^\vee)$ and $\cH(G_{n'}^{\prime\vee},\fs^{\prime\vee})$ attached to Bernstein series of enhanced $L$-parameters of $G_n$ and $G'_{n'}$, respectively.  We studied in \cite{AuHowe} the case where $\fs=[T,\triv]_{G_n}$ (with $T$ a maximal torus in $G$ and $\triv$ the trivial character of $T$) and $n'\in\{n,n+1\}$,  via the Kazhdan-Lusztig parametrization of the simple modules of $\mathcal{H}(G_n,\mathfrak{s})$ and $\mathcal{H}(G'_{n'},\mathfrak{s'})$ (in this case, we have $\fs'=[T',\triv]_{G'_{n'}}$ where $T'$ is a maximal torus in $G'$). 

\smallskip

We thank the anonymous referee for helpful suggestions that allowed to improve the readability of the paper. 

\section{Representations of finite reductive groups} \label{sec:finite}
Let $k$ be  a finite field and let $\overline{k}$ be a fixed algebraic closure of $k$. Let $\bbbG$ be a connected reductive algebraic group defined over $k$, and let $\bbG=\bbbG(k)$  denote the group of the $k$-rational points of $\bbbG$. The usual scalar product  $(\;,\;)_{\bbG}$  on the space of class functions on $\bbG$ is defined by
\begin{equation} \label{eqn:scalar product}
 (f_1,f_2)_{\bbG} := |\bbG|^{-1}\sum_{g\in \bbG}f_1(g)\,\overline{f_2(g)} .
\end{equation}
Let $\Irr(\bbG)$ be the set of equivalence classes of irreducible representations of $\bbG$. For $\pi\in\Irr(\bbG)$, we will also write $\pi$ for its character.
By Corollary~7.7 of \cite{DL}, for any $\pi\in\Irr(\bbG)$, there exists a $k$-rational maximal torus $\bbbT$ of $\bbbG$ and a character $\theta$ of $\bbT$ such that $\pi$ occurs in the Deligne-Lusztig virtual character $R_\bbbT^\bbbG(\theta)$, \textit{i.e.,} $(\pi, R_\bbbT^\bbbG(\theta))_{\bbG}\ne 0$.

\smallskip

Deligne and Lusztig introduced in \cite{DL} the following two regularity conditions:  a character $\theta$ of $\bbT$ is said to be
\begin{itemize}
\item  \textit{in general position} if its stabilizer in $(\Nor_\bbbG(\bbbT)/\bbbT)(k)$ is trivial,
\item \textit{non-singular} if it is not orthogonal to any coroot.
\end{itemize} 
If the center of $\bbbG$ is connected, then $\theta$ is non-singular if and only if it is in general position (see Proposition~5.16 in \cite{DL}).

We denote by $\bbbG^\vee$ the reductive connected group with root datum dual to that of $\bbbG$.
The $\bbG$-conjugacy classes of pairs $(\bbbT,\theta)$ as above are in one-to-one correspondence with the $\bbG^\vee$-conjugacy classes of pairs $(\bbbT^\vee,s)$ where  $\bbbT^\vee$ is a $k$-rational maximal torus of $\bbbG^\vee$ and $s$ is a semisimple element of $\bbG^\vee$ belonging to  $\bbbT^\vee$.

\smallskip

The \textit{Lusztig series} $\cE(\bbG,s)$ is defined to be the set
\[ \left\{\pi\in\Irr(\bbG)\colon
  \text{$(\pi, R_\bbbT^\bbbG(\theta))_{\bbG}\ne 0$, where $(\bbbT,\theta)_{\bbG}$ corresponds to $(\bbbT^\vee,s)_{\bbG^\vee}$}\right\}.
\]
The set $\Irr(\bbG)$  decomposes into a disjoint union
\begin{equation} \label{eqn:Lusztig-series}
\Irr(\bbG)=\bigsqcup_{(s)}\cE(\bbG,s),
\end{equation}
where $(s)$ runs over  the $\bbG^{\vee}$-conjugacy classes of semisimple elements of $\bbG^{\vee}$.  If $\theta$ is the trivial character of $\bbT$, then the representation $\pi$ is called \textit{unipotent}. 
By definition, $\cE(\bbG,1)$ is the subset of unipotent representations in $\Irr(\bbG)$. 
We denote by $\bbbG^\vee(s)$  the centralizer of $s$ in $\bbbG^\vee$, and by $\bbG^\vee(s)$ the group of $k$-rational points of $\bbbG^\vee(s)$.
\begin{remark} \label{rem:centralizer}  If $\pi\in\cE(\bbG, s)$ and $\bbbG^\vee(s)$
is a torus, then  $\pi= \pm R_\bbbT^\bbbG(\theta)$ with $\theta$ in general position (see Lemma 3.6.9 in \cite{AX-G2}).
\end{remark}
If the center of $\bbbG^\vee$ is connected then the group $\bbbG^\vee(s)$ is always connected. However,  in general it may be disconnected, and we extend the notion of Deligne-Lusztig character as follows. We denote by $\bbbG^\vee(s)^{\circ}$ the identity component of $\bbbG^\vee(s)$, and  by $\bbG^\vee(s)^\circ$ the group of $k$-rational points of $\bbbG^\vee(s)^{\circ}$. For $\bbbT^\vee$ a $k$-rational maximal torus of $\bbbG^\vee(s)^{\circ}$ and $\theta^\vee$ a character of $\bbT^\vee$, we set
\begin{equation}
R_{\bbbT^\vee}^{\bbbG^\vee(s)}(\theta^\vee):=\Ind_{\bbG^\vee(s)^{\circ}}^{\bbG^\vee(s)}(R_{\bbbT^{\vee}}^{\bbbG^\vee(s)^{\circ}}(\theta^\vee)).
\end{equation}
We define $\mathcal{E}(\bbG^\vee(s),1)$ to be the set of irreducible constituents of $R_{\bbbT^\vee}^{\bbbG^\vee(s)}(1)$.

By  Theorem~4.23 of \cite{Lubook} when the center of $\bbbG^\vee$ is connected, and by \S12 of \cite{Lus-disconnected-center} in general, there is a bijection
\begin{equation}\label{Lus-unipotent-decomposition}
    \mathcal{E}(\bbG,s)\xrightarrow{1-1} \mathcal{E}((\bbG^\vee(s),1), \quad
    \pi\mapsto\pi^\unip .
\end{equation}
The bijection \eqref{Lus-unipotent-decomposition} satisfies the following properties:
\begin{enumerate}
    \item[{\rm (1)}] It sends a Deligne-Lusztig character $R_\bbbT^\bbbG(\theta)$ in $\bbG$ (up to a sign) to a Deligne-Lusztig character $R_{\bbbT^{\vee}}^{\bbbG^\vee(s)}(1)$, where $1$ denotes the trivial character of $\bbbT^{\vee}$ (see \S12 of \cite{Lus-disconnected-center}). 
    \item[{\rm (2)}] It preserves cuspidality in the following sense:  if $s\in\bbG^\vee$, and $\mathcal{E}(\bbG,s)$ contains a cuspidal representation
$\pi_\cc$, then
\begin{itemize}
    \item[{\rm (a)}]
the largest $k$-split torus in the center of $\bbbG^\vee(s)$ coincides with the largest $k$-split torus in the center of $\bbbG^\vee$ (see (8.4.5) of \cite{Lubook}),
    \item[{\rm (b)}] the unipotent representation $\pi_\cc^\unip$ is cuspidal.
 \end{itemize}
    \item[{\rm (3)}]
The dimension of every representation $\pi\in\mathcal{E}(\bbG,s)$ is given by
\begin{equation} \label{eqn:dim tau}
\dim(\pi)=\frac{|\bbG|_{p'}}{|\bbbG^\vee(s)(k)|_{p'}}\,\dim(\pi^\unip), 
\end{equation}
where $|\bbG|_{p'}$ is the largest prime-to-$p$ factor of the order of $\bbG$ (see Remark~13.24 of \cite{DM1}).
\end{enumerate}

\subsection{Centralizers of semisimple elements in  classical groups} \label{subsec:centralizers}
Let $\bbbG$ be  a classical group of rank $n$ over $\overline{k}$, \textit{i.e.,} $\bbbG$ is a symplectic group, a special orthogonal group or an orthogonal group. Let $\bbbT_{n}\simeq\overline{k}^\times\times\cdots\times\overline{k}^\times$ be a $k$-rational maximal torus of $\bbbG$. For $s = (\lambda_1,\ldots,\lambda_{n})\in \bbT_{n}$,  we denote by $\nu_\lambda(s)$ the number of the $\lambda_i$'s which are equal to $\lambda$, and  by $\langle\lambda\rangle$ the set of all roots in $\overline{k}$ of the irreducible polynomial of $\lambda$ over $k$.  By \S1.B of \cite{AMR}, the centralizer $\Cent_\bbbG(s)$ of $s$ in $\bbbG$ decomposes as a product
 \begin{equation} \label{eqn:centr}
 \Cent_\bbbG(s)=\prod_{\langle\lambda\rangle\subset\{\lambda_1,\ldots,\lambda_{n}\}}\bbbG_{(\lambda)}(s),
 \end{equation}
where $\bbbG_{(\lambda)}(s)$ is a reductive quasi-simple group of rank equal to $\nu_\lambda(s)\cdot |\langle\lambda\rangle|$. The followings hold:
\begin{itemize}
\item[(i)] if $\bbbG=\Sp_{2n}$, then $\bbbG_{1}(s)\simeq\Sp_{2\nu_1(s)}$ and $\bbbG_{-1}(s)\simeq\Sp_{2\nu_{-1}(s)}$,
\item[(ii)] if $\bbbG=\SO_{2n+1}$, then $\bbbG_{1}(s)\simeq \SO_{2\nu_1(s)+1}$  and $\bbbG_{-1}(s)\simeq\rO^\pm_{2\nu_{-1}(s)}$,
\item[(iii)] If $\bbbG=\rO^\pm_{2n}$, then $\bbbG_{1}(s)\simeq \rO^\pm_{2\nu_1(s)}$ and $\bbbG_{-1}(s)\simeq \rO^\pm_{2\nu_{-1}(s)}$,
\item[(iv)] If $\bbbG=\rU_n$ or if $\lambda\ne\pm 1$, then $\bbbG_{(\lambda)}(s)$ is either a general linear group or a unitary group.
\end{itemize}
We set
\begin{equation}
\bbbG_{\ne}(s):=\prod_{\langle\lambda\rangle\subset\{\lambda_1,\ldots,\lambda_{m}\},\lambda\ne\pm 1}\bbbG_{(\lambda)}(s).
\end{equation}
By \eqref{eqn:centr}, we have
\begin{equation}
\Cent_{\bbbG^\vee}(s)=\bbbG_{1}(s)^\vee\times \bbbG_{-1}(s)\times\bbbG_{\ne}(s),
\end{equation} 
and hence a bijection
\begin{equation}
\begin{matrix}\cE(\Cent_{\bbG^\vee}(s),1)&\simeq &\cE(\bbG_{1}(s)^\vee,1)\times \cE(\bbG_{-1}(s),1)\times\cE(\bbG_{\ne}(s),1)\cr
\pi^\unip&\mapsto&\pi^\unip_{1}\otimes\pi^\unip_{-1}\otimes\pi^\unip_{\ne}
\end{matrix}.
\end{equation}
We obtain a one-to-one correspondence
\begin{equation} \label{eqn:Ls}
\begin{matrix}\cL_s\colon&\cE(\bbG,s)&\to&\cE(\bbG_{1}(s)^\vee,1)\times \cE(\bbG_{-1}(s),1)\times\cE(\bbG_{\ne}(s),1)\cr
&\pi&\mapsto&\pi^\unip_{1}\otimes\pi^\unip_{-1}\otimes\pi^\unip_{\ne}
\end{matrix}.
\end{equation}

\subsection{The Howe correspondence over finite fields} \label{subsec:finiteHowe}
We suppose that the characteristic of $k$ is odd.  Let $N$ be a positive integer, and let $W$ be a vector space over $k$ of dimension $2N$, equipped with a nondegenerate alternated bilinear form $\langle\,,\,\rangle$. 
A pair of reductive subgroups of $\Sp(W)=\Sp_{2N}(k)$, where each one is the centralizer of the other, is called  \textit{reductive dual pair}. We study irreducible dual pairs, because these are the building blocks of all the others.  Such pairs are of the following kinds:  
\begin{itemize}
\item[$\circ$] pairs of type I:
\begin{itemize}
\item $(\Sp_{2n}(k),\rO_{N'}(k))$ with $nN'=N$;
\item $(\rU_n(k),\rU_{n'}(k))$ with $nn'=2N$;
\end{itemize}
\item[$\circ$]
pairs of type II: 
\begin{itemize}
\item
$(\GL_n(k),\GL_{n'}(k))$ with $nn'=2N$.
\end{itemize}
\end{itemize}

A reductive dual pair $(\bbG,\bbG')$ is said to be in the \textit{stable range} (with $\bbG'$ smaller) if the defining vector space for $\bbG$ has a totally isotropic subspace of dimension greater or equal than the dimension the defining vector space for $\bbG'$, e.g. the pairs $(\Sp_{2n}(k),\rO_{2n'}(k))$ such that $n\ge 2n'$.

\smallskip

Howe introduced  a correspondence $\Theta_{n,n'}\colon \fR(\bbG_n)\to\fR(\bbG_{n'})$  between the categories of complex representations of these subgroups. It is obtained from a particular representation of $\Sp_{2N}(k)$, called the \textit{Weil} (or  \textit{oscillator})  \textit{representation}.
In order to define this  representation, we must introduce the Heisenberg group. This is the group with underlying set $H(W)= W\times k$ and product
\begin{equation}
(w_1,t_1)\cdot (w_2,t_2)=\left(w_1+w_2,t_1+t_2 + \frac{1}{2}\langle w_1,w_2\rangle\right).
\end{equation}
Let $\varrho$ be an irreducible representation of $H(W)$. Its restriction to the center $\rZ_{H(W)}\simeq k$ of $H(W)$ equals $\psi_\varrho$, for a certain character $\psi_\varrho$ of $k$.

By the Stone-von-Neumann Theorem (see Theorem 2.I.2 in \cite{MVW}), for any non-trivial character $\psi$ of $\rZ_{H(W)}$ there exists (up to equivalence) a unique irreducible representation $\varrho$ of $H(W)$ such that $\psi_\varrho=\psi$. This representation is known as the Heisenberg representation. It depends on $\psi$, so we denote it by $\varrho_\psi$.

The natural action of $\Sp(W)$ on $H(W)$ fixes the elements of its center. Hence, for a fixed character $\psi$ of $k$, the representations $\varrho_\psi$ and $g\cdot \varrho_\psi$ agree on $\rZ_{H(W)}$, for any $g\in\Sp(W)$.  The Stone-von-Neumann Theorem implies that there is an operator $\omega_\psi(g)$ verifying
\begin{equation}
\varrho_\psi(g\cdot w,t) = \omega_\psi(g)\varrho_\psi(w,t)\omega_\psi(g)^{-1}.
\end{equation}
This defines a projective representation $\omega_\psi$ of $\Sp(W)$, which can be lifted to an actual representation of $\Sp(W )$, known as the Weil representation.

The restriction of $\omega_\psi$ of $\Sp_{2N}(k)$ to $\bbG_n\times \bbG'_{n'}$ is
\begin{equation}\omega_{\bbG_n,\bbG_{n'}'}=\sum_{\pi\in\Irr(\bbG)\atop\pi'\in\Irr(\bbG')}\mult_{\pi,\pi'}\,\pi\otimes\pi',\;\;\text{where $\mult_{\pi,\pi'}\in \ZZ_{\ge 0}$.}
\end{equation}
Define $\Theta_{n'}\colon \ZZ\,\Irr(\bbG_n)\to \ZZ\,\Irr(\bbG'_{n'})$ by
\begin{equation}\Theta_{n'}(\pi):=\left\{\pi'\in\Irr(\bbG'_{n'})\,:\,\mult_{\pi,\pi'}\ne 0\right\},\quad \text{for $\pi\in\Irr(\bbG_n)$.}
\end{equation}
The occurrence of a  irreducible representation $\pi$ of $\bbG_n$ in the Howe correspondence for $(\bbG_n,\bbG'_{n'_\pi})$ 
with $n'_\pi$ minimal (i.e., such that $\Theta_{\bbG_{n'}'}(\pi)=0$ for any $n'<n'_\pi$) is referred to as the \textit{first occurrence}.

\subsubsection{The Howe correspondence for unipotent representations} \label{subsubsec:Howe-finite-unip}

Between members of a dual pair, the only ones having cuspidal unipotent representations are:
$\GL_1(k)$, $\Sp_{2(a^2+a}(k)$, $\rU_{(a^2+a)/2}(k)$ (which have a unique such representation, say $\tau_a$), and $\rO_{2a^2}(k)$ (which has two: $\tau_a^\rI$ and $\tau_a^{\rII}=\tau_a^\rI\otimes\sign$). 

\smallskip

From now on, we will only consider pairs formed by a symplectic group $\Sp_{2n}(k)$ and an orthogonal group $\rO_{2n'}^\epsilon(k)$, where $\epsilon\in\{\pm\}$, with $\rO_{2n'}^+(k)$ split and $\rO_{2n'}^-(k)$ nonsplit. The following results were established in \cite{Adams-Moy}:
\begin{enumerate}
\item if $\pi$ is a cuspidal irreducible  representation of $\bbG_n=\Sp_{2n}(k)$, then $\Theta_{n'_\pi}(\pi)$ is a singleton $\{\pi'\}$ with $\pi'$ cuspidal irreducible,
\item
if $\pi\in\Irr(\bbG_n)$ is unipotent then any $\pi'\in \Theta_{n'}(\pi)$ is unipotent,
\item the representation $\tau_a$ of $\Sp_{2(a^2+a)}(k)$ corresponds to $\tau_a^{\rII}$ if $\epsilon$ is the sign of $(-1)^a$ and to $\tau_{a+1}^{\rI}$ otherwise.
\end{enumerate}
Thus, the Howe  correspondence between cuspidal unipotent representations is describe by the function 
$\theta\colon\NN\to \NN$, defined by 
\begin{equation}
\theta(a):=
\begin{cases}
a&\text{if $\epsilon$ is the sign of $(-1)^a$}\cr
a+1&\text{otherwise.}
\end{cases}
\end{equation}
By \cite{AMR}, the Howe correspondence for unipotent representations of the dual pair $(\Sp_{2n}(k),\rO^\epsilon_{2n'}(k))$  induces a correspondence between the parabolically induced representations 
$\ii^{\Sp_{2n}(k)}_{\Sp_{2(a^2+a)}(k)\otimes \bbT}(\tau_a\otimes1)$ and 
$\ii_{\rO^\epsilon_{2\theta(a)^2}(k)\otimes \bbT'}^{\rO_{2n'}^\epsilon(k)}(\tau'_{\theta(a)}\otimes1)$, where $\tau_{\theta(a)}'\in\{\tau_{\theta(a)}^\rI,\tau_{\theta(a)}^\rII\}$, and $\bbT$, $\bbT'$ are products of $\GL_1(k)$'s. 

It induces a correspondence between the endomorphism algebras of these parabolically induced representations, and hence a correspondence  $\Omega_{N_a,N'_a}$ between irreducible representations of  Weyl groups of types $\rB_{N_a}$ and $\rB_{N'_a}$, where 
%\begin{block}{} We denote by $\omega_{n,n'}^{\unip}$ the unipotent part of $\omega_{\bbG,\bbG'}$:
%\[\omega_{n,n'}^{\unip}:=\sum_{\pi\in\Irr^{\unip}(G)\atop\pi'\in\Irr^{\unip}(\bbG')}\mult_{\pi,\pi'}\,\pi\otimes\pi'.\]
%\end{block}
$N_a:=n-a^2-a$ and $N'_a:=n-\theta(a)^2$.

The following conjectural explicit description of $\Omega_{N_a,N'_a}$ for pairs $(\Sp_{2n}(k),\rO^\epsilon_{2n'}(k))$ was formulated in \cite{AMR}, and confirmed in 
\cite{Pan3}, and in \cite{MQZ}, independently. 

Let us introduce first some definitions and notation:
$\lambda:=(\lambda_1\ge\lambda_2\ge\cdots\ge\lambda_l)$ where $\lambda_i\in\ZZ_{\ge 0}$ is called a partition of $n$ if $|\lambda|:=\lambda_1+\lambda_2+\cdots+\lambda_l=n$. We write $\lambda \dashv n$.
We denote by $\cP(n)$ set of partitions of $n$, and we
% and $\cP:=\bigcup_{n\in\NN}\cP(n)$.
write $\lambda\cup\mu$  for the partition of $n+m$ with parts  $\lambda_1,\ldots,\lambda_l,\mu_1,\ldots,\mu_l$.
The usual order on partitions is defined by
\[\lambda\le\lambda'\quad\text{if and only if}\quad \lambda_1+\cdots+\lambda_i\le \lambda_1'+\cdots+\lambda_i',\text{ for all $i\in\NN$.}\]
We define another order on partitions as follows: for $\lambda$, $\lambda'$ partitions of possibly different integers, we write
\begin{equation} \label{eqn:order}
\lambda\preceq\lambda'\quad\text{if and only if}\quad \lambda'_{i+1}\le\lambda_i\le\lambda_i' ,\text{ for all $i\in\NN$.}
\end{equation}
It says that $\lambda\preceq\lambda'$ if the Young diagram of $\lambda$ is contained in the one of $\lambda'$ and that we can go from the first to the second by adding at most one box per column. For instance, the partitions $\lambda= (4,1,1)$ of $6$, and $\lambda'= (4,4,1,1)$ of $10$ verify $\lambda\preceq\lambda'$. We observe that if two partitions $\lambda=(\lambda_1\ge\lambda_2\ge\cdots\ge\lambda_l)$  and $\lambda'=(\lambda_1\ge\lambda_2\ge\cdots\ge\lambda_{l'})$  satisfy $\lambda\preceq\lambda'$ and $|\lambda|=|\lambda'|$, then we have $\lambda=\lambda'$. Indeed,  by \eqref{eqn:order}, we have  $l\le l'$, and 
\[|\lambda|=\sum_{i=1}^l\lambda_i\leq \sum_{i=1}^l\lambda_i'=|\lambda|-(\lambda'_{l+1}+\cdots+\lambda'_{l'}),\]
thus,  $l'=l$ and $\lambda_i=\lambda'_i$ for all $i\in\{1,\ldots, l\}$.

\smallskip

Irreducible characters of a Weyl group of type $\rB_{n}$ or $\rC_{n}$ are known to be parametrized by bipartitions of $n$ (see for instance \cite{Car}). We denote by $\chi_{\xi,\eta}$ the irreducible character which  corresponds to the bipartition $(\xi,\eta)$ of $n$. 

In the rest of this section we consider pairs $(\Sp_{2n}(k),\rO^\epsilon_{2n'}(k))$ in the stable range (\textit{i.e.,} such that  that $n\ge 2n'$). We write $N_a(\zeta):=N_a-|\zeta|$ if $\zeta$ is a partition such that $|\zeta|\le N_a$. Then we have:
\begin{enumerate}
\item Cases $(\Sp_{2n},\rO^+_{2n'})$ with $a$ even and $(\Sp_{2n},\rO^-_{2n'})$ with $a$ odd: 
\[\Omega_{N_a,N'_a}=\sum_{r=0}^{\min(N_a,N'_a)}\sum_{(\xi,\zeta)\in\cP_2(r)}\sum_{\eta,\eta'}\chi_{\xi,\eta}\otimes\chi_{\xi,\eta'},\]
where the third sum is over the partitions $\eta\dashv N_a(\xi)$ and $\eta'\dashv N_a'(\xi)$ such that $\zeta\preceq\eta$ and $\zeta\preceq\eta'$.
\item Cases $(\Sp_{2n},\rO^+_{2n'})$ with $a$ odd and $(\Sp_{2n},\rO^+_{2n'})$ with $a$ even: 
\[\Omega_{N_a,N'_a}=\sum_{r=0}^{\min(N_a,N'_a)}\sum_{(\xi,\zeta)\in\cP_2(r)}\sum_{\xi',\eta'}\chi_{\xi',\eta}\otimes\chi_{\xi,\eta'},\]
where the third sum is over the partitions $\xi'\dashv N_a(\eta)$ and $\eta'\dashv N_a'(\xi)$ such that $\xi\preceq\xi'$ and $\eta\preceq\eta'$.
\end{enumerate} 

\smallskip

In general, there exist representations $\pi\in\Irr(\bbG)$ such that $\Theta_{n'}(\pi)$ contains more than one element. Hence, one may wonder if it would be possible to extract a \textit{one-to-one correspondence}.
Several approaches to this questions were considered:
\begin{itemize}
\item definition of the eta correspondence  for dual pairs $(\Sp_{2n}(k),\rO_{N'}(k))$ in the stable range in \cite{GH}, 
\item construction of a one-to-one correspondence for unipotent representations of pairs of type II  and of pairs in the stable range of the form
$(\rO_{2a^2}^\epsilon(k),\Sp_{2(a^2+a+N)}(k))$ or 
$(\Sp_{2(a^2+a+2)}(k),\rO_{2a^2+N}^\epsilon(k))$ \cite{AKP},
\item construction of a one-to-one theta correspondence for unipotent representations of irreducible pairs of type I in the stable range in \cite{Epequin},
\item extension of both eta and theta correspondences to all irreducible pairs of type I in \cite{Pan5}.
\end{itemize}

Let $(\xi',\eta')\in\cP_2(N'_a)$. We set
\begin{equation}
\Theta_{\xi',\eta'}:=\left\{\chi_{\xi,\eta}\,:\,\text{$\chi_{\xi,\eta}\otimes\chi_{\xi',\eta'}$ occurs in $\Omega_{N_a,N'_a}$}\right\}.
\end{equation}
Recall that a partition is called \textit{symplectic} if each odd part appears with even multiplicity. There is a bijection between symplectic partitions of $2n$ and unipotent conjugacy classes of $\bbbG_n=\Sp_{2n}(\overline{k})$. We denote by $\cO_\lambda$ the unipotent class associated to the symplectic partition $\lambda$.
The Springer correspondence maps $\chi_{\xi,\eta}$ to a pair $(\cO_\lambda,\rho)$ for some irreducible character $\rho$ of $\Cent_{\bbbG_n}(u)/\Cent_{\bbbG_n}(u)^\circ$, where $u\in\cO_\lambda$. In this situation we write $\lambda_{\xi,\eta}:=\lambda$. 
\begin{definition}
We say that $\chi_{\xi_\maxi,\eta_\maxi}\in \Theta_{\xi',\eta'}$ is \textit{maximal}  if 
\[\lambda_{\xi,\eta}\le \lambda_{\xi_\maxi,\eta_\maxi}\quad\text{for all $(\eta,\xi)\in\Theta_{\xi',\eta'}$.}\]
\end{definition}

\begin{remark} Since the order is not total, it is not clear a priori that a maximal representation exists, and if so, that it is unique.
\end{remark}
\begin{definition} We say that  $\chi_{\xi_\mini,\eta_\mini}\in \Theta_{\xi',\eta'}$ is \textit{minimal}  if 
\[\lambda_{\xi_\mini,\eta_\mini}\le \lambda_{\xi,\eta}\quad\text{for all $(\eta,\xi)\in\Theta_{\xi',\eta'}$.}\]
A representation is said to be \textit{extremal} if it is either maximal or minimal.
\end{definition}
Let $(\xi',\eta')\in\cP_2(N'_a)$. By \cite{Epequin},
\begin{enumerate}
\item there exists a unique maximal representation $\chi_{\xi_\maxi,\eta_\maxi}\in \Theta_{\xi',\eta'}$, it is given by
\[\xi_\maxi:=\xi'\quad\text{and}\quad \eta_\maxi:=(N_a-N_a'+\eta_1'+\eta_2',\eta_3',\ldots,\eta'_l).\]
\item there exists a unique minimal representation $\chi_{\xi_\mini,\eta_\mini}\in \Theta_{\xi',\eta'}$, it is given by
\[\xi_\mini:=\xi' \quad\text{and}\quad\eta_\mini:=(N_a-N_a')\cup\eta'.\]
\end{enumerate}

\section{Smooth representations of $p$-adic reductive groups} \label{sec:Bernstein}
Let $F$ be a non Archimedean local field and $W_F$ its abolute Weil group. We denote by $I_F,P_F\subset W_F$ the inertia group and the wild inertia group of $F$, respectively, and by $k_F$ the residual field of $F$ (a finite field with $q$ elements, where $q$ is a power of a prime number $p$). We denote by  $W_F^\tame:=W_F\rtimes I_F$ the tame Weil group of $F$.

Let $\bG$ be a connected reductive algebraic group defined over $F$. We denote by $G$ the group of the $F$-rational points of $\bG$. 
Let $L$ be a Levi subgroup of a parabolic subgroup $P$ of $G$ and let $\fX_{\nr}(L)$ denote the group of unramified characters of $L$.  
Let $\sigma$ be an irreducible supercuspidal smooth representation of $L$ and $\cO$ the set of equivalence classes of representations $L$ of the form $\sigma\otimes \chi$, with $\chi\in X_\nr(L)$. We write $\fs:=(L,\cO)_G=[L,\sigma]_G$ for the $G$-conjugacy class of the pair $(L,\cO)$ and $\fB(G)$ for the set of such classes $\fs$. We set $\fs_L:=(L,\cO)_L$.

We denote by $\fR^\fs(G)$ the full subcategory of $\fR(G)$ whose objects are  the representations $(\pi,V)$ such that every irreducible $G$-subquotient of $\pi$ has its supercuspidal support in  $\fs$. The categories $\fR^\fs(G)$ are indecomposable and split the full smooth category $\fR(G)$ in a direct product:
\begin{equation} \label{eqn:Bernstein}
\fR(G)=\prod_{\fs\in\fB(G)}\fR^\fs(G).
\end{equation}
Let $\Irr^\fs(G)$ denote the set of irreducible objects of the category $\fR^\fs(G)$. As a direct consequence of \eqref{eqn:Bernstein}, we have
\begin{equation} \label{eqn:Bersnstein series}
\Irr(G)=\prod_{\fs\in\fB(G)}\Irr^\fs(G).
\end{equation}

\subsection{Depth-zero supercuspidal representations of $p$-adic groups} \label{sec:sc}
We denote by $\cB(\bG,F)$ the Bruhat-Tits building of $G$ and by $\cB_\reduit(\bG,F):=\cB(\bG/\rZ_{\bG},F)$ the reduced  Bruhat-Tits building of $G$. We have $\cB(\bG,F)=\cB_\reduit(\bG,F)\times (X_*(\rZ_{G})\otimes_\ZZ\RR)$, where $X_*(\rZ_G)$ is the set of $F$-algebraic cocharacters of $\rZ_G$
For each point $x$ in $\cB(\bG,F)$,  we denote by $G_{x,0}$ the parahoric subgroup of $G$ associated to $x$, by $G_{x,0^+}$ the pro-$p$ unipotent radical of $G_{x,0}$, and $\bbG_{x,0}$ the quotient $G_{x,0}/G_{x,0+}$, (the points of) a connected reductive group over the residue field of $F$. Let $[x]$  denote the image of $x$ in $\cB_\reduit(\bG,F)$ and  $G_{[x]}$ the stabilizer of $[x$] under the action $G$ on $\cB_\reduit(\bG,F)$. 

A \textit{depth-zero} representation of $G$ is a representation which admits non-zero invariant vectors under the action of the pro-$p$ unipotent radical $G_{x,0+}$ of a parahoric subgroup $G_{x,0}$ of $G$. 

Let $\pi$ be any irreducible depth-zero supercuspidal representation of $G$. There exists a vertex $x\in\cB_\reduit(\bG,F)$ and an irreducible cuspidal representation $\pi_x$ of $\bbG_x$, such that the restriction of $\pi$ to $G_{x,0}$ contains the inflation $\infl(\pi_x)$ of $\pi_x$ (see \S1-2 of \cite{Morris-ENS} or  Proposition~6.6 of \cite{Moy-PrasadII}). The normalizer $\rN_{G}(G_{x,0})$ of $G_{x,0}$ in $G$ is a  totally disconnected group that is compact mod center, which, by the proof of (5.2.8) in \cite{BTII}, coincides with the fixator $G_{[x]}$ of $[x]$ under the action of $G$ on $\cB_\reduit(\bG,F)$. Let $\bpi_x$ denote an extension of $\infl(\pi_x)$ to $G_{[x]}$. Then $\pi$ is compactly induced from a representation of $G_{[x]}$: 
\begin{equation} \label{eqn:depth zero supercuspidal}
\pi=\cInd_{G_{[x]}}^G(\bpi_x).
\end{equation}
The representation $\pi$ is called \textit{unipotent} if $\pi_x$ is unipotent

\subsection{Depth-zero supercuspidal representations of $p$-adic classical groups} \label{subsec:unip-class}
\begin{example} \label{ex:unitary-padic}
Let $E$ be an unramified quadratic extension of $F$, and let $V$ be a vector space over $E$ of dimension $n$. There are two classes of hermitian forms on $V$. They are distinguished by the parity of the valuation of the determinant of the form: we denote the form by $\langle\,,\,\rangle^+$ if the valuation is even, and by $\langle\,,\,\rangle^-$ otherwise. We denote the corresponding unitary groups by $\rU^\pm_n(F)$. The group  $\rU^+_n(F)$ is quasi-split.  If $n$ is odd, then the groups $\rU^+_n(F)$ and $\rU^-_n(F)$ are isomorphic, if $n$ is even, then $\rU^-_n(F)$ is an inner form of $\rU^+_n(F)$.

The  standard maximal parahoric subgroups of $G=\rU_n^\pm(F)$ are the groups $G_{x,0}$ such that $\bbG_{x,0}$ is isomorphic to the product of two unitary groups $\rU_{n_1}(k_F)$ and $\rU_{n_2}(k_F)$ with $n_1+n_2=n$. By \S\ref{subsec:finiteHowe} and \eqref{eqn:depth zero supercuspidal}, the supercuspidal unipotent representations  of $\rU_n(F)$ are the representations $\pi=\cInd_{G_{[x]}}^L(\bpi_x)$, where $\pi_x=\tau_{a_1}\otimes\tau_{a_2}$ with $a_1^2+a_1=2n_1$ and $a_2^2+a_2=2n_2$. For $i\in\{1,2\}$, we denote by $\lambda_i$ the partition of $n_i$ defined by $\lambda_i:=(k_i, k_i-1,k_i-2,\ldots, 1,)$.
\end{example}

\begin{example} \label{ex:SOodd}
Let $G=\SO_{2n+1}(F)$.
The  standard maximal parahoric subgroups $G_{x,0}$ of $G$ have type $\rD_{n_1}\times\rB_{n_2}$ with $n_1+n_2=n$. We have 
$\bbG_{x,0}\simeq\SO_{2n_1}(k_F)\times\SO_{2n_2+1}(k_F)$. The corresponding depth-zero supercuspidal representation of $G$ in compactly 
induced from the the inflation to $G_{x,0}$ of an irreducible cuspidal representation $\pi_x^1\otimes\pi_x^2$ of $\bbG_{x,0}$. We have $\pi_x^1\in\cE(\SO_{2n_1},(s_1))$ and $\pi_x^2\in\cE(\SO_{2n_2+1},(s_2))$. By \ref{subsec:centralizers}, the group $W_{s_1}$ is of one of the following types
\begin{enumerate}
\item $\rD_{a_1^2}\times\rD_{b_1^2}$ with  $a_1^2+b_1^2=n_1$ and $a_1,b_1\in\ZZ_{\ge 0}$, 
\item $\rA^2_{a_1-1}\times\rA^2_{\ell_1-1}$ with  $a_1=\frac{c_1+c_1^2}{2}$, $\ell_1=\frac{d_1+d_1^2}{2}$, $a_1+\ell_1=n_1$ and $c_1,d_1\in\ZZ_{\ge 0}$,
\item $\rD_{a_1^2}\times\rA_{a_1-1}$ with  $a_1=\frac{c_1+c_1^2}{2}$, $a_1^2+a_1=n_1$ and $a_1,c_1\in\ZZ_{\ge 0}$,
\end{enumerate}
and $W_{s_2}$ is of one of the following types
\begin{enumerate}
\item $\rD_{a_2^2}\times\rC_{b_2^2+b_2}$ with $a_2^2+b_2^2+b_2=n_2$ and $a_2,b_2\in\ZZ_{\ge 0}$, 
\item $\rA^2_{a_2-1}\times\rA^2_{\ell_3-1}$ with  $a_2=\frac{c_2+c_2^2}{2}$, $\ell_2=\frac{d_2+d_2^2}{2}$, $a_2+\ell_1=n_2$ and $c_2,d_2\in\ZZ_{\ge 0}$,
\item $\rD_{a_2^2}\times\rA_{a_2-1}$ with  $a_2=\frac{c_2+c_2^2}{2}$, $a_2^2+a_2=n_2$ and $a_2,c_2\in\ZZ_{\ge 0}$,
\item $\rC_{b_2^2+b_2}\times\rA_{a_2-1}$ with  $a_2=\frac{c_2+c_2^2}{2}$, $b_2^2=b_2+a_2=n_2$ and $b_2,c_2\in\ZZ_{\ge 0}$.
\end{enumerate}
\end{example}
%We consider the bijections
%\begin{equation} \label{bij:SOodd1}
%(a_1^2,b_1^2)\mapsto\left(\frac{(a_1+b_1)^2}{2},\frac{(a_1-b_1)^2}{2}\right)
%\end{equation}\begin{equation} \label{bij:SOodd2}
%(a_2^2,b_2^2+b_2)\mapsto\left(\frac{(a_2+b_2)(a_2+b_2+1)}{2},\frac{(a_2-b_2)(a_2-b_2-1)}{2}\right)
%\end{equation}

\section{Affine Hecke algebras} \label{sec:AHH}
Let $\cR:=(X,R,Y,R^\vee)$ be a root datum. We denote by $W$ the group generated by the $s_\alpha$ for $\alpha \in R$. It is a finite Weyl group.  We write $S := \{s_\alpha\, : \,\alpha\in\Delta\}$. Then $(W,S)$ is a finite Coxeter system. 

We choose, for every $s \in S$, a complex number $q_s$, such that
\begin{equation} q_s = q_{s'} \quad\text{ if $s$ and $s'$ are conjugate in $W$.}
\end{equation}
Let $\bq\colon S \to\CC$ be the function $s\mapsto q_s$. We define a new $\CC$-algebra $\cH(W,\bq)$ which has a vector space basis $\{T_w \,:\, w \in W \}$. Here $T_1$ is the unit element and the $T_W$ satisfy the following quadratic relations
and braid relations 
\begin{equation} 
(T_s+1)(T_s-q_s)=0 \quad\text{and}\quad 
\underbrace{T_sT_{s'}T_s\cdots}_{\text{$m(s,s')$ terms}} = \underbrace{T_{s'}T_sT_{s'}\cdots}_{\text{$m(s,s')$ terms}} \quad\text{for any $s,s'\in S$.}
\end{equation}
Fix $q\in\RR_{>1}$  and let $\lambda,\lambda^*\colon R\to\CC$ be functions such that
\begin{itemize}
\item 
 if $s_\alpha$ and $s_\beta$ are conjugate in $W$, then $\lambda(\alpha)=\lambda(\beta)$ and $\lambda^*(\alpha)=\lambda^*(\beta)$, 
 \item if $\alpha^\vee\notin 2Y$, then $\lambda^*(\alpha)=\lambda(\alpha)$.
 \end{itemize}
 We note that $\alpha^\vee\in 2Y$ is only possible for short roots $\alpha$ in a type $\rB$ component of the root system $R$. 
 
For $\alpha\in R$ we write
\begin{equation}
q_{s_\alpha} :=q^{\lambda(\alpha)} \quad
\text{and (if $\alpha^\vee\in m_{\max}^\vee$)} \;\;
q_{s_\alpha'}:=q^{\lambda^*(\alpha)}.
\end{equation}
We denote by $\ell$ the usual length function on $W$.
Let $\cH(W,\bq)$ denote the Iwahori-Hecke algebra of $W$. It has a basis $\{T_w : w \in W \}$ such that
\[
\begin{array}{cccl}
T_w T_v & = & T_{wv} & \text{if } \ell (w) + \ell (v) = \ell (wv), \\
(T_{s_\alpha} + 1) (T_{s_\alpha} - q_{s_\alpha}) & 
= & 0 & \text{if } \alpha \in\Delta.
\end{array}
\]
Let $\{\theta_x\,:\,x\in X\}$ denote the standard basis of $\CC[X]$. 
Then the affine Hecke algebra $\cH(\cR,\lambda,\lambda^*,q)$ is the vector space 
$\CC [X] \otimes_\CC \cH(W,\bq)$ such that $\CC [X]$ and $\cH(W,\bq)$ are embedded as subalgebras, and
for $\alpha\in\Delta$ and $x\in X$:
\[
\theta_x T_{s_\alpha} - T_{s_\alpha} \theta_{s_\alpha (x)}
=\left( (q^{\lambda (\alpha)} - 1) + \theta_{-\alpha} 
\left(q^{\frac{\lambda (\alpha)+\lambda^* (\alpha)}{2}}-q^{\frac{\lambda (\alpha)-\lambda^* (\alpha)}{2}}\right)\right) \frac{\theta_x - 
\theta_{s_\alpha (x)}}{\theta_0 - \theta_{\! -2\alpha}}.\]
 When $ \alpha^\vee \notin 2 Y$, the cross relation simplifies to
\[\theta_x T_{s_\alpha} - T_{s_\alpha} \theta_{s_\alpha (x)} =(q^{\lambda (\alpha)} -1) \frac{\theta_x - 
\theta_{s_\alpha (x)}}{\theta_0 - \theta_{-\alpha}}.\]

\subsection{Affine Hecke algebras and Bernstein blocks} \label{subsec:AHA-B}
Let $\fs=[L,\sigma]_G\in\fB(G)$. 
%If $\Pi^\fs$ is a progenerator $\fR^\fs(G)$, then there is a natural equivalence from $\fR^\fs(G)$ to the category $\End_G(\Pi^\fs)$-Mod of right $\End_G(\Pi^\fs)$-modules, namely $\cV\mapsto\Hom_G(\Pi^\fs,\cV)$.
We write $\fs_L:=[L,\sigma]_L$ and denote by $W^\fs$ the extended finite Weyl group $\Nor_G(\fs_L)/L$.
The set of roots of $G$ with respect to $L$ contains a root system $\Sigma^\fs$, namely the set of roots for which the associated Harish-Chandra $\mu$-function has a zero on the orbit $\fX_\nr(L)\cdot\sigma$. This induces a semi-direct decomposition
\begin{equation}
W^\fs= W(\Sigma^\fs)\rtimes\fR^\fs,
\end{equation}
where $W(\Sigma^\fs)$ is the Weyl group of $\Sigma^\fs$ and $\fR^\fs$ is the stabilizer in $W^\fs$ of the set of positive roots.
We have
\begin{equation}
L^1=\bigcap_{\chi\in\fX_\nr(L)}\ker(\chi),
%\quad\text{and}\quad L^2_\sigma:=\bigcap_{\chi\in\fX_\nr(L,\sigma)}\ker(\chi).
\end{equation}
and we suppose that the restriction of $\sigma$ to $L^1$ is multiplicity free.
%, and we set $\cR:=(\Sigma^\fs^\vee,L^2_\sigma/L^1,(L^2_\sigma/L^1)^\vee,\Sigma^\fs)$.
By \cite{Solleveld-endomorphism-algebra}, there exists a  $2$-cocycle 
\begin{equation}
\natural^\fs\colon \fR^{\fs}\times \fR^{\fs}\longrightarrow\CC[\fX_\nr(L)\cdot\sigma]^\times,
\end{equation}
such that we have a bijection
\begin{equation} \label{eqn:Sol}
\Irr^\fs(G) \longrightarrow \Irr(\cH(G,\fs)),
\end{equation}
where $\cH(G,\fs)$ is the twisted extended affine Hecke algebra of the form
\begin{equation}
\cH(G,\fs):=\cH(\cR,\lambda,\lambda^*,q)\rtimes\CC[\fR^\fs,\natural^\fs].
\end{equation}

\subsection{Affine Hecke algebras and enhanced $L$-parameters} \label{subsec:AHA-Lp}
We denote by $W'_F:=W_F\times\SL_2(\CC)$  the Weil-Deligne group of $F$ and by $G^\vee$ the complex connected reductive group with root datum dual to that of $G$. The results in \cite{AMS1}, \cite{AMS3} apply to an arbitrary group $G$, but, for the simplicity of the exposition, we suppose here that $G$ is a pure inner twist of an $F$-split group $G^*$.

A \textit{Langlands parameter} -- or $L$-parameter -- is  then a morphism  $\varphi\colon W_F'\to G^\vee$ such that
\begin{itemize}
\item $\varphi|_{\SL_2(\CC)}$ is morphism of algebraic groups,
\item $\varphi(w)$ is a semisimple element of $G^\vee$, for any $w\in W_F$.
\end{itemize}
Let $\rZ_{G^\vee}$ denote the center of $G^\vee$. 

For $\varphi$ a given  $L$-parameter, we define
\begin{equation} \label{eqn:Sphi}
S_\varphi:=
\Cent_{G^\vee}(\varphi(W'_F)).
\end{equation}

An \textit{enhanced $L$-parameter} is a pair $(\varphi,\rho)$ where 
$\varphi$ is an  $L$-parameter for $G$ and $\rho\in\Irr(\cS_\varphi)$, with $\cS_\varphi:=S_\varphi/S_\varphi^\circ$.
For $\varphi$ a given  $L$-parameter,  the representation $\rho$ is called an \textit{enhancement} of $\varphi$. We define
an action of $G^\vee$ on the set of enhanced $L$-parameters by:
\begin{equation}
g\cdot (\varphi,\rho):=(g\varphi g^{-1},{}^g\rho),\;\;\text{for $g\in G^\vee$,}
\end{equation}
where ${}^g\rho\colon h\mapsto \rho(g^{-1}hg)$.
We denote by $\Phi_{\enh}$ set of $G^\vee$-conjugacy classes of enhanced $L$-parameters, and by $\Phi_{\enh}(G)$ the subset formed by the ones that are relevant for $G$.

By applying \eqref{eqn:Sphi} to the restriction $\varphi|_{W_F}$ of $\varphi$ to $W_F$,  we define the (possibly disconnected) complex reductive group
\begin{equation} \label{eqn:Gphi}
\cG_\varphi:=\Cent_{G^\vee}(\varphi(W_F\times\{1\})).
\end{equation}
We denote by $\cG_\varphi^\circ$ the identity component of $\cG_\varphi$, and we set
\begin{itemize}
\item 
$u=u_\varphi:=\varphi\left(1,\left(\begin{smallmatrix} 1&1\cr 0&1
\end{smallmatrix}\right)\right)$ (a unipotent element);
\item $A_{\cG_\varphi}(u_\varphi) := \pi_0 (\Cent_{\cG_\varphi} (u))$. We have $\cS_\varphi\simeq A_{\cG_\varphi}(u_\varphi)$.
\end{itemize}

An enhanced $L$-parameter $(\varphi,\rho)\in\Phi_\enh$ is called \textit{cuspidal} if the following properties hold:
\begin{itemize}
\item
$\varphi$ is discrete (i.e., $\varphi(W_F')$ is not contained in any proper Levi subgroup of $G^\vee$), 
\item
$(u_\varphi,\rho)$ is a \textit{cuspidal pair} in $\cG_\varphi$ (see \cite{LuIC}, and Remark~4.6, Definition~4.7 and Definition~4.11 in \cite{AuPune}).
\end{itemize}
We denote by $\Phi_{\enh,\cusp}(G)$ the set of $G^\vee$-conjugacy of cuspidal enhanced $L$-parameters for $G$.

\begin{example}
Let $G=\Sp_{2n}(F)$ and let $\iota_{G^\vee}\colon G^\vee\hookrightarrow \GL_{2n+1}(\CC)$ be the standard embedding. 
For every $L$-parameter $\varphi$ for $G$, we define 
\begin{equation}
I_\varphi:=\left\{\tau\in \Irr(W_F)\,:\,\text{$\tau$ occurs in $\iota_{G^\vee}\circ\varphi|_{W_F}$}\right\}.
\end{equation}
We have $I_\varphi=I^\rO_\varphi\sqcup I^\rS_\varphi \sqcup I^\GL_\varphi$, where $I^\rO_\varphi$ (resp. $I^\rS_\varphi$) is the subset of $I_\varphi$ formed by the representations that are orthogonal (resp. symplectic), and $I_\varphi^\GL$ is the maximal subset of $I_\varphi$ which is formed by representations $\tau$ that are not selfdual  and satisfy $\tau^\vee\in I_\varphi^\GL$ for every $\tau\in I_\varphi^\GL$.
Hence, we have the following decomposition
\begin{equation}
\iota_{G^\vee}\circ\varphi=\bigoplus_{\tau\in I_\varphi^\rO}\tau\boxtimes M_\tau\oplus \bigoplus_{\tau\in I_\varphi^\rS}\tau\boxtimes M_\tau\oplus \bigoplus_{\tau\in I_\varphi^\GL}(\tau\oplus\tau^\vee)\boxtimes M_\tau
\end{equation}
where $M_\tau$ is  a multiplicity space of $\tau$. Let $m_\tau$ denote the dimension of $M_\tau$. We have
\begin{equation}
\cG_\varphi\simeq\prod_{\tau\in I_\varphi^\rO} \Sp_{m_\tau}(\CC)\times\prod_{\tau\in I_\varphi^\rS} \rO_{m_\tau}(\CC)\times\prod_{\tau\in I_\varphi^\GL} \GL_{m_\tau}(\CC).
\end{equation}
By \cite{Mou}, if the enhanced $L$-parameter $(\varphi,\rho)$  is cuspidal, then we have
\begin{equation}
m_\tau=\begin{cases}
d_\tau^2+ d_\tau \;\;\;\text{(with $d_\tau\in \ZZ_{\ge 0}$)}&\text{if $\tau\in I_\varphi^\rO$}\cr
d_\tau^2 \;\;\;\text{(with $d_\tau\in \ZZ_{>0}$)}&\text{if $\tau\in I_\varphi^\rS$}\cr
1&\text{if $\tau\in I_\varphi^\GL$}.
\end{cases}
\end{equation}
\end{example}

Let $\cC\in\Unip(\cG_\varphi^\circ)$ and let $\cE$ irreduciblebe an  $\cG_\varphi^\circ$-equivariant local system on $\cC$.
By \cite{LuIC}, the $\IC$-sheaf  $\cF_\rho:=\IC(\cC,\cE_\rho)$ occurs as a summand of $\ri_{\cL\subset\cP}^{\cG^\circ}(\IC(\cC_{\cusp},\cE_{\cusp}))$, for some triple $(\cP,\cL,(\cC_{\cusp},\cE_{\cusp}))$, where $\cP$ is a parabolic subgroup of $\cG^\circ$ with Levi subgroup $\cL$ and $(\cC_{\cusp},\cE_{\cusp})$ is a cuspidal unipotent pair in $\cL$. 
Moreover,  the triple $(\cP,\cL,\cC_{\cusp},\cE_{\cusp})$ is unique up to $\cG^\circ$-conjugation. 

Let $\rho^\circ\in\Irr(A_{\cG^\circ}(u))$.  The cuspidal support of  $(u,\rho^\circ)$, denoted by $\Sc^{\cG^\circ}(u,\rho^\circ)$,  is defined to be
\begin{equation}
(\cL,(v,\rho^\circ_{\cusp}))_{\cG^\circ},\quad\text{where $v\in\cC_\cusp$ and $\rho^\circ_\cusp\leftrightarrow \cE_{\cusp}$.}\end{equation}
We set $\cT:=\rZ_{\cL}^\circ$ and $\cM:=\Cent_{\cG}(\cT)$. The cuspidal support of  $(u,\rho)$ is a (well-defined) triple $(\cM,v,\rho_{\cusp})_\cG$, where $\rho^\circ_{\cusp}$ occurs in the restriction of $\rho_{\cusp}$ to $A_{\cG^\circ}(u)$.
The cuspidal support of $(\varphi,\rho)$ is defined to be
\begin{equation} \label{eqn:Sc Galois side1}
\Sc(\varphi,\rho):=(\Cent_{{}^LG}(\cT),(\varphi_v,\rho_{\cusp})).
\end{equation}

Recall from Section 3.3.1 of \cite{Hai} that the group of unramified characters of $L$ is naturally
isomorphic to $(\rZ_{L^\vee} \rtimes I_F )^\circ_{W_F}$. We consider this as an
object on the Galois side of the local Langlands correspondence and we write
\begin{equation}
\fX_\nr ({}^L L) := (\rZ_{L^\vee \rtimes I_F})_{W_F}^\circ .
\end{equation}
Given $(\varphi,\rho) \in \Phi_{\enh}(L)$ and $z \in (\rZ _{L^\vee \rtimes I_F})_{W_F}$,
we define $(z \cdot\varphi,\rho) \in \Phi_{\enh}(L)$ by
\begin{equation}
z \cdot\varphi:= \varphi \text{ on } I_F \times \SL_2 (\CC) \text{ and }
(z \cdot\varphi)(\Fr_F) := \tilde z \varphi(\Fr_F) ,
\end{equation}
where $\tilde z \in (\rZ _{L^\vee \rtimes I_F})_{W_F}$ represents $z$.

An inertial equivalence class for $\Phi_{\enh}(G)$ is the $G^\vee$-conjugacy class $\fs^\vee$ of a pair $({}^LL,\fs^\vee_L)$, where $L$
is a Levi subgroup of $G$, ${}^LL=L^\vee\rtimes W_F$ and $\fs^\vee_L$ is a $\fX_\nr ({}^L L)$-orbit in $\Phi_{\enh,\cusp} (L)$. 

The \textit{Bernstein series} $\Phi_\enh^{\fs^\vee}(G)$ associated to $\fs^\vee\in\fB(G^\vee)$ is defined to be the fiber of $\fs^\vee$ under the map $\Sc$ defined in \eqref{eqn:Sc Galois side1}. By \cite{AMS1}, the set $\Phi_\enh(G)$ of $G^\vee$-conjugacy classes of enhanced $L$-parameters for $G$ is partitioned,  analogously to \eqref{eqn:Bersnstein series}, as
\begin{equation} \label{eqn:decPhi_e}
\Phi_\enh(G)=\prod_{\fs^\vee\in\fB(G^\vee)}\Phi^{\fs^\vee}(G).
\end{equation}
In \cite{AMS3}, we canonically associated an affine Hecke algebra (possibly extended with a
finite $R$-group) $\cH(G^\vee,\fs^\vee)$ to  every Bernstein series $\Phi_\enh^{\fs^\vee}(G)$ of enhanced Langlands
parameters for $G$. While we considered only the first case of \eqref{eqn:Sphi}, our construction applies equally well to both the other cases. We showed that the simple modules of this algebra are naturally in
bijection with the elements of the Bernstein series $\Phi_{\enh}^{\fs^\vee}(G)$ and that the set of central characters of
the algebra is naturally in bijection with the collection of cuspidal supports of these
enhanced Langlands parameters. 

We summerise the construction of $\cH(G^\vee,\fs^\vee)$. By applying \eqref{eqn:Sphi} to the restriction $\varphi|_{I_F}$ of $\varphi$ to the inertia group $I_F\subset W_F$,  we define the (possibly disconnected) complex reductive group
\begin{equation} \cJ_\varphi:=\Cent_{G^\vee}(\varphi(I_F\times\{1\})).\end{equation}
Let $R(\cJ^\circ,\cT)$ be the set of $\alpha\in X^*(\cT)\backslash\{0\}$ which appear in the adjoint action of $\cT$ on the Lie algebra of $\cJ_\varphi^\circ$.
By Proposition 3.9 in \cite{AMS3}, $R(\cJ^\circ,\cT)$ is a root system. We denote by $W_{\fs^\vee}^\circ$ its Weyl group.
Let $W_{\fs^\vee}:=\Nor_{G^\vee}(\fs^\vee)/L^\vee$. We have $W_{\fs^\vee}=W_{\fs^\vee}^\circ\rtimes\Fr_{\fs^\vee}$, where 
\[\Fr_{\fs^\vee}:=\left\{w\in W_{\fs^\vee}\,:\,w (R(\cJ^\circ,\cT)^+)\subset R(\cJ^\circ,\cT)^+\right\}.\]
We define a root datum
\[\cR_{\fs^\vee}:=(R_{\fs^\vee},X^*(T_{\fs^\vee}),R_{\fs^\vee}^\vee,X_*(T_{\fs^\vee})),\]
where $T_{\fs^\vee}\simeq \fs^\vee_L$ and 
\[R_{\fs^\vee}=\left\{m_\alpha \,\alpha\,:\,\alpha\in R(\cJ^\circ,\cT)_{\mathrm{red}}\subset X^*(T_{\fs^\vee})\right\},\]
with $m_\alpha\in \ZZ_{>0}$. The group $W_{\fs^\vee}$ acts on $\cR_{\fs^\vee}$.
In \cite{AMS3} we defined $W_{\fs^\vee}$-invariant functions
\begin{equation}\lambda\colon m_{\fs^\vee}\to \QQ_{>0}\quad\text{and}\quad 
\lambda^*\colon\{m_\alpha \alpha\in m_{\fs^\vee}:m_\alpha \alpha\in 2X_*(T_{\fs^\vee})\}\to \QQ.
\end{equation}
We choose semisimple subgroups $\cJ_i \subset \cJ^\circ$, normalized by $\rN_{G}(\fs)$,
such that the derived group $\cJ^\circ_\der$ is the almost direct product of the $\cJ_i$,
and the multiplication map
\begin{equation}\label{eq:0.1}
m_{\cJ^\circ} \colon\rZ_{\cJ^\circ}^\circ \times \cJ_1 \times \cdots \times \cJ_d \to \cJ^\circ
\end{equation}
is a surjective group homomorphism with finite central kernel.
It induces  a decomposition
\begin{equation}\label{eq:0.6}
\Lie(\cJ) = \rZ_{\Lie(\cJ ^\circ)} \oplus \Lie(\cJ_1) \oplus \cdots \oplus \Lie(\cJ_d).\end{equation}
We obtain an orthogonal, $W_{\fs^\vee}$-stable decomposition
\begin{equation}\label{eq:0.2}
R(\cJ^\circ,\cT) = m_1 \sqcup \cdots \sqcup m_d, \quad\text{where $m_i:=R(\cJ_i,\cJ_i\cap \cT)$.}
\end{equation}
We let $\vec r:= (m_1, \ldots, m_d)$ be an array of variables, corresponding to \eqref{eq:0.1} in the sense that $m_i$ is relevant for $\cJ_i$ and $m_i$ only.  We have
 \begin{equation}
\cH(G^\vee,\fs^\vee)=\cH(\cR,\lambda,\lambda^*,\vec r)\rtimes\CC[\fR^{\fs^\vee},\natural^{\fs^\vee}],
\end{equation}
where $\cH(\cR,\lambda,\lambda^*,\vec r)$ is a generalized version (with $d$ indeterminates) of the algebras considered in \S\ref{sec:AHH} and $\natural^{\fs^\vee}\colon \fR^{\fs^\vee}\times\fR^{\fs^\vee}\to\CC^\times$ is $2$-cocycle.
The cocycle is trivial when $G$ is an inner twist of $\GL_n(F)$, a pure inner twist. of a quasi-split classical group or of the group $\GSpin_n(F)$ (see \cite{AMS3}, \cite{AMS4}), and when $G$ is the exceptional group of type $\rG_2$ (see \cite{AX-Hecke}).

\section{Howe correspondence} \label{sec:Howe}
\subsection{A correspondence between Hecke algebras of $p$-adic groups} \label{subsec:Howe-Hecke}
We will consider here only the case of dual pairs formed by a $p$-adic symplectic group and a $p$-adic orthogonal group of even dimension. The other cases can be treated in a similar way.
Let $n$ and $n'$ be two fixed non-negative integers, and let $W_n$ be a symplectic space of (even) dimension $2n$ over $F$. The corresponding group of isometries  is the symplectic group $G_n:=\Sp(W_n)$.
Let $V_{n'}$ be a quadratic space of dimension $2n'$ over $F$ (i.e., a space endowed with a non-degenerate symmetric
$F$-bilinear form). We denote by $\eta_{n'}$ the character of $F^\times$  associated to $F(\sqrt\Delta_{n'})/F$, where $\Delta_{n'}$ is the discriminant of the form. The orthogonal group $G_{n'}':=\rO(V_{n'})$ is the corresponding group of
isometries, and $(G_n,G_{n'}')$ is a dual pair in $\Sp(W_{2nn'})$. Fixing a non-trivial additive character
$\psi$ of $F$, we obtain the so-called Weil representation $\omega$ of the metaplectic cover $\Mp(W_{2nn'})$ of $\Sp_{W_{2nn'}}$. We can define a splitting $G _n\times G'_{n'}\to\Mp(W_{2nn'})$ (see \cite{Kudla2}). By means of this splitting we obtain a Weil representation of $G_n\times  G'_{n'}$. which we denote $\omega_{n,n'}$.  

For every $\pi\in\Irr(G_n)$,  the maximal $\pi$-isotypic quotient of $\omega_{n,n'}$ is of the form $\pi\otimes\Theta_{n'}(\pi)$, where $\Theta_{n'}(\pi)$ is a smooth representation of $G'_{n'}$, called the \textit{full theta lift} of $\pi$ to $W_{n'}$. The representation $\Theta_{n'}(\pi)$, when non-zero, has a unique irreducible quotient, denoted by
$\theta_{n'}(\pi)$, which is called the \textit{small theta lift} of $\pi$. Similarly, if $\pi'\in\Irr(G'_{n'})$, the maximal $\pi'$-isotypic quotient
of $\omega_{n,n'}$ is of the form $\Theta_n(\pi')\otimes\pi'$, where $\Theta_n(\pi')$ is a smooth representation of $G_n$, called the full theta lift of $\pi'$ to $V_n$. The
representation $\Theta_n(\pi')$, when non-zero, has a unique irreducible quotient, denoted $\theta_n(\pi')$, which is called the small theta lift of $\pi'$.
These statements were first formulated by Howe in \cite{How}, and were proved by Waldspurger in \cite{Walds} when the residual characteristic of $F$ is odd, and by Gan and Takeda \cite{GT} in general.

We have $V_{n'}=V^\an\oplus H^{n'}$, where $V^\an$ is an anisotropic $F$-vector subspace, and $H$ is the hyperbolic plane. For $m'\ge 0$, let
$V_{m'}:=V^\an\oplus H^{m'}$. The collections $\cbT:=\left\{W_{m}\,:\, m\ge 0\right\}$  and $\cbT':=\left\{V_{m'} \,:\, m' \ge 0\right\}$ are called Witt tower of vector spaces. 
One can then consider a tower of the theta correspondence associated to the tower of reductive dual pairs $(G_n, G'_{m'})_{m'\ge 0}$. For an irreducible smooth representation $\pi$ of $G_n$, we thus have the representation $\theta_{m'}(\pi)$. The smallest non-negative integer $n'_{\min}$ such that $\theta_{n'_{\min}}(\pi)\ne 0$ is called the \textit{first occurrence index} of $\pi$ for the Witt tower $\cbT'$, and the representation  $\theta^0(\pi):=\theta_{n'_{\min}}(\pi)$ is called the first occurrence of $\pi$ for this Witt tower. Similarly, one can then consider a tower of the theta correspondence associated to the tower of reductive dual pairs $(G_{m}, G'_{n'})_{r\ge 0}$. For an irreducible smooth representation $\pi'$ of $G'_n$, we have the representation $\theta_{m'}(\pi)$. The smallest non-negative integer $n_{\min}'$ such that $\theta_{n_{\min}}(\pi')\ne 0$ is called the \textit{first occurrence index} of $\pi'$ for the Witt tower $\cbT$, and the representation  $\theta^0(\pi'):=\theta_{n_{\min}}(\pi)$ is called the first occurrence of $\pi'$ for this Witt tower. 
By Chapter~3 in \cite{MVW}, such $n_{\min}'$ and $n_{\min}$ exist and $n_{\min}'\le n'$ and $n_{\min}\le n$. Moreover, we have $\Theta_{m'}(\pi)\ne 0$ for any $m'\ge n'_{\min}$ and $\Theta_{m}(\pi')\ne 0$ for any $m\ge n_{\min}$.
By Theorem~2.5 in \cite{Kudla} (see also Chapter~3 in \cite{MVW}), if the representation $\pi$ is supercuspidal then $\Theta_{n'_{\min}}(\pi)$ is  irreducible (and thus is equal to $\theta_{n'_{\min}}(\pi)$) and supercuspidal. Similarly, if $\pi'$ is supercuspidal then $\Theta_{n_{\min}}(\pi')$ is  irreducible (and thus is equal to $\theta_{n_{\min}}(\pi')$) and supercuspidal.

A Levi subgroup $L_n$ of a parabolic subgroup $P_n$ of $G_n$ is isomorphic to $G_{n_0}\times \GL_{n_1}(F)\times\cdots\times\GL_{n_r}(F)$, where $n_0+n_1+\cdots+n_r=n$, and any $\pi\in\Irr(G_n)$ is isomorphic to a subquotient of the parabolic induced representation $\ii_{L_n,P_n}(\sigma)$ where $\sigma\in\Irr_\scc(L_n)$. We have $\sigma\simeq\pi_{n'_0}\otimes\tau_1\otimes\cdots\otimes\tau_r$, where $\pi_{n'_0}\in\Irr_\scc(G'_{n'_0})$ and $\tau_i\in\Irr_\scc(\GL_{m_i}(F))$ for $1\le i\le r$.
A Levi subgroup $L_{n'}'$ of a parabolic subgroup $P_{n'}'$ of $G_{n'}'$ is isomorphic to $G'_{n'_0}\times \GL_{n'_1}(F)\times\cdots\times\GL_{n'_{r'}}(F)$,  where $n'_0+n'_1+\cdots+n'_{r'}=n'$ and any $\pi'\in\Irr(G_n')$ is isomorphic to a subquotient of the representation $\ii_{L_n',P_n'}(\sigma')$, where $\sigma'\simeq\pi_{n'_0}'\otimes\tau_1'\otimes\cdots\otimes\tau_{r'}'\in\Irr_\scc(L'_{n'})$ such that $\pi_{n'_0}'\in\Irr_\scc(G_{n'_0}')$ and $\tau_j'\in\Irr_\scc(\GL_{n'_j}(F))$ for $1\le j\le r'$.  The $G_n$-conjugacy class $(L_n,\sigma)_{G_n}$  of the pair $(L_n,\sigma)$ is uniquely determined by $\pi$, is called the supercuspidal support of $\pi$, and denoted by $\Sc(\pi)$.  Similarly, the $G_{n'}'$-conjugacy class of the pair $(L_{n'}',\sigma_{n'}')$ is uniquely determined by $\pi'$,  is called the supercuspidal support of $\pi'$,  and denoted by $\Sc(\pi')$. We set 
\begin{equation}
\overline n:=n-n_0\quad\overline n':=n'-n'_0,\quad d:=n-n',
\quad\text{and}\quad d_0:=n_0-n'_0.
\end{equation} 
By \cite{Kudla},  if $\Theta_{n'}(\pi)\ne \{0\}$ and $\pi':=\theta_{n'}(\pi)$, then 
\begin{enumerate}
\item if $\overline n'\le \overline n$, the supercuspidal support of $\pi'$ is
\[\left(L'_{n'},\pi_{n'_0}'\otimes\eta_{n'}^{-1}\tau_1\otimes\cdots\otimes\eta_{n'}^{-1}\tau_{\ell}\otimes\nu_F^{d-1}\otimes\nu_F^{d-2}\otimes\cdots\otimes\nu_F^{d_0}\right)_{G'_{n'}},
\]
where $\pi'_{n'_0}:=\theta^0(\pi_{n_0})$.
\item if $\overline n'\ge \overline n$, the supercuspidal support of $\pi$ is
\[\left(L_n,\pi_{n_0}\otimes\eta_{n'}\tau_1'\otimes\cdots\otimes\eta_{n'}\tau_{r'}'\otimes\eta_{n'}\nu_F^{-d-1}\otimes\eta_{n'}\nu_F^{-d-2}\otimes\cdots\otimes\eta_{n'}\nu_F^{-d_0}\right)_{G_n}.
\]
where $\pi_{n_0}:=\theta^0(\pi'_{n'_0})$.
\end{enumerate}

Let $L\in\{L_n,L'_{n'}\}$. 
Recall that a character of $L$ is said to be unramified if it is trivial on all the compact subgroups of $L$. We denote by $\fX_\nr(L)$ the group of unramified characters of $L$. 
\begin{proposition} \label{prop:series}
For any $\pi\in\Irr^\fs(G_n)$, where $\fs=[L_n,\sigma]_{G_n}$, we have $\theta_{n'}(\pi)\in\Irr^{\theta(\fs)}(G_n)$, where 
\begin{equation}
\theta(\fs):=[L_{n'}',\sigma']_{G_{n'}'}.
\end{equation}
\end{proposition}
\begin{proof} Let $\pi\in\Irr^\fs(G_n)$, where $\fs=[L_n,\sigma]_{G_n}$ and $\pi':=\theta_{n'}(\pi)$. Thus, we have $\Sc(\theta_{n'}(\pi))=(L'_{n'},\sigma')_{G'_{n'}}$. We set $\fs':=[L'_{n'},\sigma']_{G'_{n'}}$.
Let $\chi\in \fX_\nr(L_n)$, and let $\pi_\chi\in\Irr(G_n)$ such that $\Sc(\pi_\chi)=(L_n,\chi\otimes\sigma)_{G_n}$. The character $\chi$ has the following form:  $\chi=\triv\otimes\chi_1\otimes\cdots\otimes\chi_r$, where $\chi_i\in\fX_\nr(\GL_{n_i}(F))$ for $1\le i\le r$.
By applying the previous formula (case (1))  to the representation $\pi_\chi$ we see that, if $\overline n'\le \overline n$, then the supercuspidal support of $\theta_{n'}(\pi_\chi)$ is
\[\begin{matrix}
\left(L'_{n'},\pi_{n'_0}'\otimes\chi_1\eta_{n'}^{-1}\tau_1\otimes\cdots\otimes\chi_r\eta_{n'}^{-1}\tau_{r}\otimes\nu_F^{n'-n}\otimes\cdots\otimes\nu_F^{\overline n'-\overline n-1}\right)_{G'_{n'}}\cr
=(L'_{n'},(\chi\otimes
\underbrace{\triv\otimes\cdots\otimes\triv}_{\textrm{$n_0-n'_0$ times}})\otimes\sigma')_{G'_{n'}}.
\end{matrix}
\]
The character $\chi\otimes\triv^{\otimes (n_0-n'_0)}$ of $L'_{n'}$  being  unramified,  we have  $\theta_{n'}(\pi_\chi)\in\Irr^{\fs'}(G'_{n'})$. Since every representation in $\Irr^\fs(G_n)$ is of the form $\pi_\chi$ for some $\chi\in \fX_\nr(L_n)$, we get
\[\theta_{n'}(\Irr^\fs(G_n))\subset\Irr^{\fs'}(G'_{n'}).\]
On the other hand, if $\overline n'\ge \overline n$, by definition of the representation $\pi_\chi$ its supercuspidal support  is  $(L_n,\chi\otimes\sigma)_{G_n}$, hence, by the formula (case (2)) applied to $\pi$, it  equals 
\[
\left(L_n,\pi_{n_0}\otimes\eta_{n'}\chi_1\tau_1'\otimes\cdots\otimes\eta_{n'}\chi_{r'}\tau_{r'}'\otimes\eta_{n'}\nu_F^{n'-n}\otimes\cdots\otimes\eta_{n'}\nu_F^{\overline m-\overline n-1}\right)_{G_n}.
\]
Thus, by applying the formula (case (2)) to $\pi_\chi$, we get
\[\Sc(\theta_{n'}(\pi_\chi))=(L'_{n'},\pi'_{n'_0}\otimes\chi_1\tau_1'\otimes\cdots\otimes\chi_{r'}\tau_{r'}')_{G'_{n'}}=(L'_{n'},\chi'\otimes\sigma')_{G'_{n'}},\]
where $\chi':=\triv\otimes\chi_1\otimes\cdots\otimes\chi_{r'}\in\fX_\nr(L'_{n'})$.
Hence, $\theta_{n'}(\pi_\chi)\in \Irr^{\fs'}(G'_{n'})$.
\end{proof}

\begin{theorem} \label{thm:corrAHA}
For every $\fs\in\fB(G)$, the Howe correspondence $\pi\mapsto\theta_{n'}(\pi)$ for the reductive dual pair $(G_n,G'_{n'})$ induces a correspondence $E\mapsto \theta(E)$ between subsets of simple modules of the extended affine Hecke algebras $\cH(G_n,\fs)$ and $\cH(G'_{n'},\theta(\fs))$.
\end{theorem}
\begin{proof}
By \cite{Heiermann}, the algebras $\cH(G_n,\fs)$ for any $\fs\in\fB(G_n)$ and $\cH(G'_{n'},\fs')$ for any $\fs'\in\fB(G'_{n'})$ are extended affine Hecke algebras (that is, the $2$-cocycles are trivial in the case of  classical groups).
Then the result follows from the combination of Proposition~\ref{prop:series} with \eqref{eqn:Bersnstein series}.
\end{proof}

\subsection{A correspondence between completions of Hecke algebras} \label{subsec:completions}
To study representations of $G$ it is often useful to consider various group algebras of $G$. There is the Hecke algebra $\cH(G)$, defined to be the convolution algebra of locally constant, compactly supported functions $f\colon G\to\CC$. The category $\fR(G)$ is equivalent to the category of of nondegenerate $\cH(G)$-modules. By letting $G$ act on $\cH(G)$ by left translation, we obtain from \eqref{eqn:Bernstein} a decomposition 
\begin{equation} \label{eqn:decH}
\cH(G)=\bigoplus_{\fs\in\fB(G)}\cH(G)^\fs,
\end{equation}
with $\cH(G)^\fs\in\fR^\fs(G)$. The spaces $\cH(G)^\fs$ are two-sided ideals of $\cH(G)$. 

From the point of view of noncommutative geometry, for the study of tempered representations of $G$, it is interesting to use the reduced $C^\ast$-algebra $C^\ast_\rr(G)$, whose spectrum coincides with the tempered dual $\Irr^\temp(G)$ of G. Analogously, we have  the following decomposition of $C_\rr^\ast(G)$:
\begin{equation}
C_\rr^\ast(G)=\bigoplus_{\fs\in\fB(G)}C_\rr^\ast(G)^\fs.
\end{equation}

\begin{theorem} \label{thm:Cstar}
We suppose that $n'=n$ or $n'=n+1$.  Then the theta correspondence $\pi\mapsto\theta_{n'}(\pi)$ for the reductive dual pair $(G_n,G'_{n'})$ induces a correspondence  between subsets of simple modules of $C_\rr^\ast(G)^\fs$ and $C_\rr^\ast(G)^{\theta(\fs)}$.
\end{theorem}
\begin{proof}
By Theorem~1.2  in \cite{GT0}, if $\pi\in\Irr^\temp(G_n)$, then $\theta_{n'}(\pi)\in\Irr^\temp(G'_{n'})$. On the other hand, by Theorem~6.2 in \cite{Solleveld-completion}, when $G$ is $G_n$ or $G'_{n'}$, then $C_\rr^\ast(G)^\fs$ is Morita equivalent to the reduced $C^\ast$-completion of $\cH(G,\fs)$ for each $\fs\in\fB(G)$. Thus, the result follows from Theorem~\ref{thm:corrAHA}.
\end{proof}

\begin{remark}
Mesland and \c{S}eng\"un already observed in \cite{MS} that certain cases of the theta correspondence can be described by using objects of $C^\ast$-algebraic nature. More precisely, they constructed Hilbert $C^\ast$-bimodules over $C^\ast$-algebras of groups in equal rank dual pairs and showed that the Rieffel induction functors implemented by these modules coincide with the theta correspondence.
\end{remark}

\subsection{The Howe correspondence for depth-zero representations} \label{subsec:depth-zero-class}
In  \cite{Pan0}, Pan proved that the theta correspondence  preserves depth zero representations
for type I reductive dual pairs $(G,G')$ over $p$-adic fields of odd residual characteristic and is compatible with the theta correspondence for finite reductive dual pairs in the following sense.
Let $(\pi,\cV)$ be a depth-zero representation of $G$. There exists a point $x$ of $\cB(\bG,F)$ such that $\pi$ has non-zero invariant vectors under $G_{x,0+}$ and an irreducible cuspidal representation $\pi_x$ of $\bbG_{x,0}$ which is contained in  $\cV^{G_{x,0+}}$. If $(\pi',\cV')\in\Irr(G')$ corresponds to $(\pi,\cV)$ via the theta correspondence, then there exists a point $x'$ of $\cB(\bG',F)$ such that $\pi'$ has non-zero invariant vectors under $G'_{x',0+}$ and an irreducible cuspidal representation $\pi_{x'}$ of $\bbG'_{x',0}$ which is contained in  $\cV^{G'_{x',0+}}$, and the representation $\pi_{x'}$  corresponds to $\pi_x$ via the theta correspondence for reductive dual pair $(\bbG_{x,0},\bbG'_{x',0})$ over $k_F$.

\subsection{The Howe correspondence for supercuspidal representations} \label{subsec:Howe-sc}

\subsubsection{The Howe correspondence for positive-depth supercuspidal representations} \label{susubbsec:positive-depth-class} Loke and Ma gave in \cite{LM} a description of the local theta correspondence between tame supercuspidal representations in terms of the supercuspidal data. Tame supercuspidal representations were constructed in \cite{Yu}. 
Let $\cD(G_{n'})$ (resp. $\cD(G'_{n'})$) be a set of supercuspidal data for $G_n$ (resp. $G'_{n'}$). For a given supercuspidal Yu datum $\Sigma$, let  $[\pi_\Sigma]$ denote the isomorphism class of the supercuspidal representation attached to $\Sigma$. In \cite{HM} an equivalence relation $\sim$ on $\cD(G_{n})$ was defined
so that the map from $\overline\cD(G_{n}):= \cD(G_{n})/\sim$ to $\Irr_\scc(G_{n})$ given by $[\Sigma] \mapsto [\pi_\Sigma]$ is a bijection. We set
\[\overline\cD_{\cbT}:=\bigcup_{n}\overline\cD(G_n).\]
Using the moment maps and theta correspondences over finite fields,  a \textit{theta lift for supercuspidal data} was defined in \cite{LM}. It is a map \begin{equation}
\vartheta_{n',\cbT'}\colon \overline\cD(G_{n'})\hookrightarrow \overline\cD_{\cbT'}.
\end{equation}
Suppose that  $\Sigma\in\cD(G_{n})$ and $[\Sigma']:=\vartheta_{n,\cbT'}([\Sigma])\in\overline\cD_{n'}(G')$ for some $n'$. By \cite{LM}, we have $\theta_n(\pi_{\Sigma})=\pi_{\Sigma'}'$, and, conversely, if $\pi$ and $\pi'$ are irreducible supercuspidal representations of $G_n$ and $G_{n'}'$, respectively, then there exists $\Sigma\in\cD(G_n)$ such that  $\pi=\pi_{\Sigma}$ and $\pi'=\pi'_{\Sigma'}$, where $[\Sigma']:=\vartheta_{n',\cbT'}([\Sigma])$.

\subsection{A correspondence between Hecke algebras for enhanced $L$-parameters} \label{subsec:Howe-eLp}
Let  $\fs=[L_n,\sigma]_{G_n}\in\fB(G_n)$. When $F$ has characteristic zero, a local Langlands correspondence $\sigma\mapsto(\varphi_\sigma,\rho_\sigma)$ is available via endoscopy, due to Arthur  and M{\oe}glin (see \cite{Art} and \cite{Moe}). By \cite{GV}, this result extends to the case where $F$ has positive characteristic. 

We set 
\begin{equation}
\fs^\vee:=[L_n^\vee,(\varphi_\sigma,\rho_\sigma)]_{G_n^\vee}.
\end{equation} 
By \cite{AMS3}, for $G\in\{G_n,G'_{n'}\}$, the algebra $\cH(G,\fs)$ is canonically isomorphic to the algebra $\cH(G^\vee,\fs^\vee)$.
Thus, by Theorem~\ref{thm:corrAHA}, the Howe correspondence $\pi\mapsto\theta_{n'}(\pi)$ for the reductive dual pair $(G_n,G'_{n'})$ induces a correspondence $E\mapsto \theta(E)$ between subsets of simple modules of the extended affine Hecke algebras $\cH(G_n^\vee,\fs^\vee)$ and $\cH(G_{n'}^{\prime\vee},\theta(\fs)^\vee)$, and hence a correspondence \begin{equation}
(\varphi,\rho)\mapsto (\theta_{n'}(\varphi),\theta_{n'}(\rho))
\end{equation}
between subsets of the series $\Phi^{\fs^\vee}_\enh(G_n)$ and $\Phi^{\theta(\fs)^\vee}_\enh(G_{n'}')$.

\end{document}